\documentclass[11pt,reqno]{amsart}

\usepackage[utf8]{inputenc}
\usepackage[T1]{fontenc}
\usepackage{lmodern}
\usepackage{amsmath,amssymb,amsthm,mathtools}
\usepackage{graphicx}
\usepackage[usenames,dvipsnames]{xcolor}
\usepackage[colorlinks=true,linkcolor=NavyBlue,citecolor=OliveGreen,urlcolor=NavyBlue]{hyperref}
\usepackage{microtype}

\graphicspath{{figs/}}

\numberwithin{equation}{section}

\newtheorem{theorem}{Theorem}[section]
\newtheorem*{theoremA}{Theorem A}
\newtheorem*{theoremB}{Theorem B}
\newtheorem*{theoremC}{Theorem C}
\newtheorem{proposition}[theorem]{Proposition}
\newtheorem{lemma}[theorem]{Lemma}
\newtheorem{corollary}[theorem]{Corollary}

\theoremstyle{definition}
\newtheorem{definition}[theorem]{Definition}

\theoremstyle{remark}
\newtheorem{remark}[theorem]{Remark}

\newcommand{\R}{\mathbb{R}}
\newcommand{\C}{\mathbb{C}}
\newcommand{\Sph}{\mathbb{S}}
\newcommand{\Sb}{\Sph^3_\tau}
\newcommand{\base}{\Sph^2\!\left(\tfrac12\right)}
\newcommand{\g}{g_\tau}
\newcommand{\met}[2]{\langle #1,#2\rangle}

\newcommand{\dd}{\,\mathrm{d}}
\DeclareMathOperator{\Iso}{Isom}
\DeclareMathOperator{\tr}{tr}

\DeclareMathOperator{\diag}{diag}

\begin{document}

\title[Solitons of the MCF in Berger spheres]{Solitons of the mean curvature flow in Berger spheres:\\ Hopf tubes and rotationally invariant surfaces}

\author{Hiuri Reis}
\address{Instituto de Matem\'atica e Estat\'istica, Universidade Federal de Goi\'as, 74690-900, Goi\^ania--GO, Brazil}
\email{hiuri\_reis@ufg.br}

\author{Adriana Cintra}
\address{Instituto de Matem\'atica e Estat\'istica, Universidade Federal de Goi\'as, 74690-900, Goi\^ania--GO, Brazil}
\email{adriana.cintra@ufg.br}

\author{Guilherme Sabo}
\address{Universidade de Bras\'ilia, Departamento de Matem\'atica, 70910-900, Bras\'ilia--DF, Brazil}
\email{guipaes92@gmail.com}


\begin{abstract}
We classify the solitons of the mean curvature flow in the Berger spheres
$\Sb=\mathbb{E}(4,\tau)$: surfaces $\Sigma\subset\Sb$ evolving by ambient isometries,
$H=\g(X,N)$ for a Killing field $X$. For the two natural classes of invariant surfaces
--- Hopf tubes and surfaces invariant under $(z,w)\mapsto(z,e^{it}w)$ --- we reduce the
soliton equation to autonomous planar ODE systems for rotating curves with density, and
show a weighted momentum-type quantity is strictly monotone along every nontrivial
profile. Together with a first-variation flux argument and existence-uniqueness for the
reduced ODEs, this gives, in each class: (0) the nontrivial complete solitons form an
explicit one-parameter family, up to ambient isometry, indexed by the unique degenerate
point of the profile; (i) the only compact solitons are minimal and trivial, the
Clifford torus being the unique nontotally geodesic invariant one; (ii) every
noncompact complete soliton is an immersed spiral converging exponentially, at both
ends, to the Clifford torus; (iii) there is exactly one exceptional soliton up to
isometry --- a \emph{yin--yang Hopf tube} through a Hopf fiber, and, in the rotational
class, an axis-crossing disk and a fiber-crossing cylinder; (iv) the ends wind around
the Clifford torus iff the rotation speed is below an explicit threshold
($\alpha_*=8$ for Hopf tubes, independent of $\tau$; $a_*=4\sqrt{1+\tau^2}/\tau$ for
rotational solitons). The Hopf-tube analysis is $\tau$-independent and lifts the
spherical curve-shortening solitons of dos Reis--Tenenblat to every Berger sphere;
numerical phase portraits and stereographic plots are provided.
\end{abstract}

\maketitle
\setcounter{tocdepth}{1}
\tableofcontents

\section{Introduction}

A one-parameter family of immersions $F_t\colon\Sigma^2\to (M^3,g)$ evolves by the
\emph{mean curvature flow} (MCF) if
\begin{equation}\label{eq:MCF}
  \partial_t F = \vec H = H N ,
\end{equation}
where $N$ is a unit normal and $H=k_1+k_2$ denotes the trace of the shape operator.
Solitons --- solutions that evolve by symmetries of the ambient space --- play a
central role in the analysis of the flow: they arise as blow-up models at
singularities, as barriers, and as the natural ``rigid'' solutions against which the
general behavior of the flow is measured. In Euclidean space this program is classical:
self-shrinkers govern type-I singularities \cite{ColdingMinicozzi2012}, translators
model type-II blow-ups, and the rotating and dilating solitons of the planar curve
shortening flow were completely classified by Halld\'orsson \cite{Halldorsson2012}.
A complementary source of explicit examples comes from symmetry reductions of the
flow itself: dos Reis and Tenenblat \cite{ReisSantosTenenblat2018} show that every
isoparametric hypersurface of a space form evolves, under MCF, through the family of
its parallel hypersurfaces, reducing the flow to an explicit ODE and yielding exact
extinction times --- a reduction in the same spirit as the equivariant one we use
below.

In a Riemannian ambient space without dilations or translations the only available
symmetries are isometries, so the natural notion is the following. Let $X$ be a Killing
field on $(M,g)$ with flow $(\Phi_t)$. An oriented surface $\Sigma\subset M$ is a
\emph{soliton with respect to $X$} if
\begin{equation}\label{eq:soliton}
  H=g(X,N)\qquad\text{along }\Sigma .
\end{equation}
Then $t\mapsto\Phi_t(\Sigma)$ is, up to tangential diffeomorphisms, a solution of
\eqref{eq:MCF} (see \cite{HungerbuehlerSmoczyk2000} and Lemma~\ref{lem:solflow} below);
$\Sigma$ is \emph{trivial} when $X$ is everywhere tangent to $\Sigma$, in which case
$\Sigma$ is a minimal surface translated along itself. Solitons of curve shortening
and mean curvature flows in non-Euclidean ambients have attracted considerable
attention in recent years; in particular the second-named line of research most
relevant to us is the classification of the soliton curves of the curve shortening
flow on the round $2$-sphere by dos Reis and Tenenblat \cite{ReisTenenblat2019} and its
hyperbolic counterpart \cite{SilvaTenenblat2021}; for the full mean curvature flow
(rather than curve shortening) the analogous rotators and translators have very
recently been classified in hyperbolic $3$-space by de Lima, Ramos and dos Santos
\cite{DeLimaRamosDosSantos2026}, who also exhibit a helicoidal family of constant
mean curvature rotators extending a result of Halld\'orsson's in $\R^3$.

In this paper we carry this program into the \emph{Berger spheres}, the simplest
homogeneous deformations of the round $3$-sphere. We take the unit sphere
$\Sph^3\subset\C^2$ and, for $\tau>0$, the metric
\begin{equation}\label{eq:bergerintro}
  \g=\met{\cdot}{\cdot}+(\tau^2-1)\,\met{\cdot}{V}\met{\cdot}{V},
  \qquad V(z,w)=(iz,iw),
\end{equation}
where $\met{\cdot}{\cdot}$ is the round metric and $V$ is the Hopf field. The space
$\Sb=(\Sph^3,\g)$ is the homogeneous space $\mathbb{E}(4,\tau)$: the Hopf map
$\pi\colon\Sb\to\base$ is a Riemannian submersion with totally geodesic fibers of
length $2\pi\tau$ and bundle curvature $\tau$, over the $2$-sphere of curvature $4$
\cite{Daniel2007}. For $\tau\neq1$ the isometry group is four-dimensional, with
identity component $\mathrm U(2)$ \cite{Torralbo2010,Torralbo2012}, and
\emph{every} Killing field is, up to conjugation, of the form
\[
  X_{a,b}(z,w)=(iaz,\,ibw),\qquad a,b\in\R .
\]
Since $\Sb$ is compact, there are no translating or shrinking solitons (contrast the
Euclidean shrinking tori of \cite{Angenent1992}); every soliton is a \emph{rotator}
with respect to some $X_{a,b}$.

Two classes of surfaces interact best with these symmetries, and they are the classes
we classify: \emph{Hopf tubes} $\pi^{-1}(\gamma)$, invariant under the fiberwise circle
action generated by $V$; and \emph{rotationally invariant surfaces}, invariant under
the circle action $(z,w)\mapsto(z,e^{it}w)$, whose fixed-point set is the great circle
$\{w=0\}$ and whose ``axis'' behaves like the rotation axis of a surface of revolution.
Both classes reduce the soliton equation \eqref{eq:soliton} to an autonomous ODE in a
two-dimensional orbit space, and we treat both reductions in the unified form of
\emph{rotating curves with density} (Lemma~\ref{lem:reduction}): the profile curve
$\gamma$ satisfies
\[
  \kappa=\hat g(Y,n)+\partial_n f ,
\]
where $\hat g$ is the orbital metric, $f$ is the logarithm of the orbit length, and
$Y$ is the projected rotation field. The key structural fact
(Lemma~\ref{lem:monotone}) is that the weighted momentum
\[
  \Phi=e^{f}\,\hat g(Y,T)
\]
is \emph{strictly increasing} along every nontrivial soliton profile. This
monotonicity, an ODE analogue of Huisken-type monotonicity adapted to rotators,
excludes all periodic behavior at once and localizes the limit sets.

\subsection*{Main results}
Throughout, ``soliton'' means a complete, connected surface satisfying
\eqref{eq:soliton} for the indicated Killing field, and $\mathcal T_{\mathrm{Cl}}
=\{|z|=|w|\}$ denotes the Clifford torus, which is minimal in every $\Sb$.

\begin{theoremA}[compact solitons; Theorem \ref{thm:compact}]
Let $\Sigma\subset\Sb$ be a compact soliton with respect to a Killing field $X$. Then
$\Sigma$ is minimal and $X$ is everywhere tangent to $\Sigma$; that is, every compact
soliton is trivial. In particular, the only compact soliton in the Hopf and rotational
classes with respect to a field $X_{a,b}$ acting nontrivially on the class
($\alpha=a-b\neq0$, resp.\ $a\neq0$) is the Clifford torus.
\end{theoremA}

\begin{theoremB}[Hopf tubes; Theorem \ref{thm:hopfclass}]
Fix $X=X_{a,b}$ with $\alpha=a-b>0$. A Hopf tube $\pi^{-1}(\gamma)$ is a soliton with
respect to $X$ if and only if $\gamma$ solves the curve shortening flow soliton
equation $\kappa=\met{Y_\alpha}{n}$ on $\base$, where $Y_\alpha$ generates the rotations
of speed $\alpha$ about the poles; the condition is independent of $\tau$. This
equation is already classified, on the round sphere, by dos Reis and Tenenblat
\cite{ReisTenenblat2019}; transported to $\base$, their theorem gives:
\begin{enumerate}
\item existence-uniqueness for the reduced ODE yields, up to the isometries of
      $\base$, a one-parameter family of nontrivial complete soliton curves, indexed
      by the colatitude of their (unique) turning point;
\item the only closed soliton curve is the equator, whose tube is the Clifford torus;
\item every other complete soliton curve is a nonclosed immersed spiral, along which
      the momentum $h=\hat g(Y_\alpha,T)$ is strictly increasing between its limits
      $\mp\alpha/2$; it has exactly one point of meridional tangency or polar crossing,
      it crosses the poles at most once, and both of its ends converge exponentially,
      with rate $\alpha/4$, to the equator; the corresponding Hopf tubes are complete
      immersed cylinders whose two ends converge to the Clifford torus;
\item the ends wind while oscillating across the equator if $\alpha<8$, and approach it
      monotonically from one side if $\alpha\geq8$;
\item up to ambient isometry there is exactly one soliton curve through a pole ---
      the spherical \emph{yin--yang curve}; its Hopf tube is the unique soliton
      Hopf tube containing a whole fiber.
\end{enumerate}
\end{theoremB}

\begin{theoremC}[rotationally invariant solitons; Theorem \ref{thm:rotclass}]
Fix $X=X_{a,b}$ with $a\neq0$ and let $\Sigma$ be invariant under
$(z,w)\mapsto(z,e^{it}w)$. The soliton equation depends only on $a$ and reduces to an
autonomous system in the orbit disk $B$, endowed with the metric
$\hat g=\dd\chi^2+q(\chi)\dd\theta^2$,
$q=\tau^2\cos^2\chi/(\cos^2\chi+\tau^2\sin^2\chi)$, and the density
$\ell=\sin\chi\,(\cos^2\chi+\tau^2\sin^2\chi)^{1/2}$. Moreover:
\begin{enumerate}
\item existence-uniqueness for the reduced ODE yields, up to the residual rotations
      and reflection of $B$, a one-parameter family of nontrivial complete profiles,
      indexed by the value $\chi_0\in[0,\pi/2]$ of their (unique) degenerate point,
      with the two exceptional solitons of item (5) as its endpoints;
\item the only closed profile, hence the only compact soliton, is the Clifford torus
      $\{\chi=\pi/4\}$;
\item along every other complete profile the weighted momentum
      $\Phi=\tfrac{a\tau}{2}\sin(2\chi)\sin\beta$ is strictly increasing between
      $\mp a\tau/2$; consequently each profile meets the axis circle, the polar fiber,
      or a radial tangency exactly once, and both of its ends converge exponentially,
      with rate $\mu/2$, $\mu=a\tau/\sqrt{1+\tau^2}$, to the Clifford torus;
\item the ends spiral around the Clifford torus if $a<a_*=4\sqrt{1+\tau^2}/\tau$ and
      approach it monotonically if $a\geq a_*$;
\item up to ambient isometry there are exactly two exceptional solitons: the
      \emph{axis-crossing soliton disk}, an immersed plane meeting the fixed-point
      circle $\{w=0\}$ orthogonally at one point, and the \emph{fiber-crossing soliton
      cylinder}, containing the fiber $\{z=0\}$; all remaining solitons are immersed
      cylinders staying in the open solid tori bounded by these configurations.
\end{enumerate}
\end{theoremC}

Theorems B and C exhibit the Clifford torus as a universal attractor for invariant
solitons of the MCF in every Berger sphere, in sharp contrast with the Euclidean
picture where soliton ends escape to infinity. Another notable feature is that the
Hopf-tube problem is completely independent of the Berger parameter $\tau$ ---
solitons of the two-dimensional base lift isometrically in an essentially topological
way --- while the rotational problem sees $\tau$ both in its ODE and in its
oscillation threshold $a_*$; as $\tau\to0$ (collapse to the base) $a_*\to\infty$ and
every rotational soliton end winds, while as $\tau\to1$ we recover the round sphere
with $a_*=4\sqrt2$.

All figures in this paper were produced from high-accuracy numerical integrations of
the systems \eqref{eq:hopfsys} and \eqref{eq:rotsys}. Full computational details
behind every formula in this paper are collected in a companion set of notes.

\subsection*{Organization} Section \ref{sec:berger} describes the Berger metric, its
isometries and Killing fields. Section \ref{sec:solitons} contains the general
results: the equivariant reduction, the monotonicity lemma, Theorem A, and, in
\S\ref{sec:tools}, an explicit statement (with references) of the four tools from the
qualitative theory of planar ODEs used throughout the classification, aimed at a
reader more familiar with geometry than with dynamical systems. Section
\ref{sec:hopf} proves Theorem B; since the Hopf-tube soliton equation coincides,
$\tau$-independently, with the curve-shortening-soliton equation on the round sphere
already classified by dos Reis and Tenenblat \cite{ReisTenenblat2019}, we cite their
theorem for the phase-plane classification rather than reproving it, and detail only
the consequences genuinely new in the three-dimensional Berger setting. Section
\ref{sec:rot} proves Theorem C for the rotational class, which was not previously
classified; there we give the reduced ODE's derivation and its phase-portrait
classification in full detail, as the model instance of the method.

\subsection*{Acknowledgments} An AI assistant (Claude, Anthropic) was used
to help check English grammar and wording in the preparation of this manuscript;
all mathematical content, computations and conclusions are the authors' own.

\section{The Berger sphere, its isometries and Killing fields}\label{sec:berger}

\subsection{The metric}
Let $\Sph^3=\{(z,w)\in\C^2 : |z|^2+|w|^2=1\}$ carry the round metric
$\met{\cdot}{\cdot}$ of curvature $1$, and let $V(z,w)=(iz,iw)$ be the \emph{Hopf
field}, the unit tangent to the fibers of the Hopf fibration. For $\tau>0$ the
\emph{Berger metric} is
\begin{equation}\label{eq:berger}
 \g(U,W)=\met{U}{W}+(\tau^{2}-1)\met{U}{V}\met{W}{V} .
\end{equation}
We write $\Sb=(\Sph^3,\g)$; for $\tau=1$ this is the round sphere, and for
$\tau\neq1$ a Berger sphere proper. All results below hold for every $\tau>0$.

The Hopf map
\begin{equation}\label{eq:hopfmap}
  \pi(z,w)=\Bigl(z\bar w,\ \tfrac12\bigl(|z|^2-|w|^2\bigr)\Bigr)\in\C\times\R,
\end{equation}
is a Riemannian submersion from $\Sb$ onto the sphere $\base$ of radius $\tfrac12$
(curvature $4$), with totally geodesic circle fibers of $\g$-length $2\pi\tau$; the
vertical unit field is $\xi=V/\tau$ and satisfies
$\overline\nabla_U\,\xi=\tau\,U\wedge\xi$ for all $U$, where $\wedge$ is the
$\g$-cross product. Thus $\Sb=\mathbb E(4,\tau)$ in the notation of \cite{Daniel2007};
its sectional curvature equals $\tau^2$ on planes containing $\xi$ and $4-3\tau^2$ on
horizontal planes (companion notes).

It is convenient to use the coordinates
\begin{equation}\label{eq:coords}
 (z,w)=\bigl(\cos\chi\, e^{i\theta},\ \sin\chi\, e^{i\psi}\bigr),
 \qquad \chi\in[0,\tfrac\pi2],\ \theta,\psi\in\R/2\pi\mathbb{Z},
\end{equation}
in which $V=\partial_\theta+\partial_\psi$ and
\begin{align}
\g={}&\dd\chi^{2}
 +\cos^{2}\chi\,\bigl(\sin^2\chi+\tau^{2}\cos^{2}\chi\bigr)\dd\theta^{2}\notag\\
 &+\sin^{2}\chi\,\bigl(\cos^2\chi+\tau^{2}\sin^{2}\chi\bigr)\dd\psi^{2}
 +2(\tau^{2}-1)\sin^{2}\chi\cos^{2}\chi\,\dd\theta\,\dd\psi .
\label{eq:metriccoords}
\end{align}

\subsection{Isometries and Killing fields}
The unitary group $\mathrm U(2)$ acts on $\C^2$ preserving both $\met{\cdot}{\cdot}$
and $V$ (because it commutes with multiplication by $i$), hence acts by isometries
on every $\Sb$. Together with the anti-holomorphic isometry $c(z,w)=(\bar z,\bar w)$,
which satisfies $dc(V)=-V\circ c$, this exhausts the isometry group when $\tau\neq1$.

\begin{proposition}[{see \cite{Torralbo2010,Torralbo2012}}]\label{prop:isometries}
For $\tau\neq1$, $\Iso(\Sb)$ has dimension $4$ and identity component
$\mathrm U(2)$; moreover $\Iso(\Sb)=\mathrm U(2)\cup c\,\mathrm U(2)$.
Every Killing field is of the form
\[
 X_A(p)=A\,p,\qquad A\in\mathfrak u(2),
\]
and a basis of the Lie algebra of Killing fields is
\begin{align*}
 V   &=X_{iI}=(iz,iw)=\partial_\theta+\partial_\psi, &
 X_1 &=X_{\diag(i,-i)}=(iz,-iw)=\partial_\theta-\partial_\psi,\\
 X_2 &=(w,-z), &
 X_3 &=(iw,iz).
\end{align*}
\end{proposition}

Since every $A\in\mathfrak u(2)$ is unitarily conjugate to a diagonal matrix, every
Killing field is conjugate in $\Iso(\Sb)$ to exactly one
\begin{equation}\label{eq:Xab}
 X_{a,b}(z,w)=(iaz,ibw)=a\,\partial_\theta+b\,\partial_\psi,
 \qquad a\ge b ,
\end{equation}
whose flow is $\Phi^{a,b}_t(z,w)=(e^{iat}z,\,e^{ibt}w)$. Hence in the study of
solitons we may, and do, restrict to the fields \eqref{eq:Xab}. Under the Hopf map,
\begin{equation}\label{eq:projection}
  \pi\circ\Phi^{a,b}_t=R_{(a-b)t}\circ\pi ,
\end{equation}
where $R_s$ is the rotation of $\base$ of angle $s$ about the poles
$(0,\pm\tfrac12)$: the projected Killing field is $Y_\alpha=\alpha\,\partial_v$ with
\emph{relative speed} $\alpha=a-b$, while the Hopf field $V=X_{1,1}$ projects to zero.

\subsection{Solitons}
\begin{definition}\label{def:soliton}
An oriented surface $\Sigma\subset\Sb$ with unit normal $N$ and mean curvature
$H=\tr S$ ($S$ the shape operator) is a \emph{soliton with respect to the Killing
field $X$} if $H=\g(X,N)$ on $\Sigma$. The soliton is \emph{trivial} if
$\g(X,N)\equiv0$, in which case $\Sigma$ is minimal and invariant under the flow
of $X$.
\end{definition}

\begin{lemma}[{\cite{HungerbuehlerSmoczyk2000}}]\label{lem:solflow}
If $\Sigma$ is a soliton with respect to $X$, then $\Sigma_t=\Phi_t(\Sigma)$ is a
solution of the mean curvature flow up to tangential reparametrization, and
conversely.
\end{lemma}

\begin{remark}\label{rem:gauge}
If $K$ is a Killing field everywhere tangent to $\Sigma$, then $X$ and $X+cK$ induce
the same soliton equation on $\Sigma$. This gauge freedom will reduce the effective
parameter space in both invariant classes: for Hopf tubes only $\alpha=a-b$ matters,
and in the rotational class only $a$.
\end{remark}

\section{Equivariant reduction and general results}\label{sec:solitons}

Both classes of surfaces studied in this paper are invariant under an isometric
circle action, and the soliton equation reduces to the orbit space in a unified way.
Let $K$ be a Killing field on $\Sb$ generating an isometric action of $\Sph^1$, let
$M^\circ$ be the open set where $K\neq0$, and let $B=M^\circ/\Sph^1$ carry the orbital
distance metric $\hat g$. Write
\[
 \ell=|K|_{\g},\qquad f=\log\ell \colon B\to\R ,
\]
for the orbit-length density. If $\Sigma$ is a $K$-invariant surface, its projection
is a curve $\gamma\subset B$, the \emph{profile}, with unit tangent $T$, unit normal
$n$ and geodesic curvature $\kappa$ in $(B,\hat g)$.

\begin{lemma}[equivariant reduction]\label{lem:reduction}
Let $\Sigma=\pi_K^{-1}(\gamma)$ be a $K$-invariant surface. Then the mean curvature of
$\Sigma$ is
\begin{equation}\label{eq:redH}
  H=\kappa-\partial_n f
\end{equation}
along $\gamma$. Moreover, if $X$ is a Killing field commuting with $K$, then
$\g(X,N)=\hat g(Y,n)$ where $Y$ is the projection of $X$ to $B$, and hence $\Sigma$ is
a soliton with respect to $X$ if and only if the profile satisfies the
\emph{rotating-curve-with-density equation}
\begin{equation}\label{eq:redsoliton}
  \kappa=\hat g(Y,n)+\partial_n f .
\end{equation}
\end{lemma}

\begin{proof}
Let $e=K/\ell$ be the unit orbit direction, and lift $T,n$ horizontally. Since $K$ is
Killing, $\overline\nabla_K K=-\tfrac12\nabla|K|^2$, so
$\g(N,\overline\nabla_e e)=-\partial_N\log\ell=-\partial_n f$; this is the principal
curvature of $\Sigma$ in the orbit direction. The horizontal part of the second
fundamental form projects to the geodesic curvature of $\gamma$ in the orbital
metric, giving \eqref{eq:redH}. For the soliton term, decompose $X$ into its
components tangent and orthogonal to $K$; the first is tangent to $\Sigma$, and the
second projects to $Y$, so $\g(X,N)=\hat g(Y,n)$ because $N$ is horizontal. Details
are carried out in the companion notes.
\end{proof}

The pair $(Y,f)$ always satisfies the compatibility condition that makes the
following momentum estimate work: $Y$ is Killing on $(B,\hat g)$ and preserves the
density, $Y(f)=0$, because $X$ preserves $|K|$.

\begin{lemma}[weighted momentum monotonicity]\label{lem:monotone}
Let $\gamma$ be a solution of \eqref{eq:redsoliton} parametrized by arclength, where
$Y$ is Killing for $\hat g$ and $Y(f)=0$. Then
\begin{equation}\label{eq:monotone}
  \frac{\dd}{\dd s}\Bigl(e^{f}\,\hat g(Y,T)\Bigr)=e^{f}\,\hat g(Y,n)^2\;\ge 0,
\end{equation}
with equality on an interval only where $\hat g(Y,n)=0$. In particular
$\Phi=e^f\hat g(Y,T)$ is strictly increasing along every profile which is not an
orbit of the flow of $Y$.
\end{lemma}

\begin{proof}
Since $Y$ is Killing, $\hat g(\nabla_T Y,T)=0$, so
$\hat g(Y,T)'=\hat g(Y,\kappa n)$. Using \eqref{eq:redsoliton},
\[
\begin{aligned}
 \Phi'&=e^f\bigl[(\partial_T f)\hat g(Y,T)+\kappa\,\hat g(Y,n)\bigr] \\
      &=e^f\hat g(Y,n)^2+e^f\bigl[(\partial_T f)\hat g(Y,T)+(\partial_n f)\hat g(Y,n)\bigr] \\
      &=e^f\hat g(Y,n)^2+e^f\,df(Y)\;=\;e^f\hat g(Y,n)^2 ,
\end{aligned}
\]
since $df(Y)=Y(f)=0$.
If $\hat g(Y,n)\equiv0$ on an interval then $\gamma$ is tangent to $Y$, hence an
integral curve of $Y/|Y|$, i.e.\ an orbit.
\end{proof}

Lemma \ref{lem:monotone} is the engine of this paper: it forbids closed and
recurrent profiles and produces sharp counting statements for the special
crossings, since $\Phi$ vanishes exactly on the degeneracy locus
$\{\ell=0\}\cup\{|Y|=0\}\cup\{T\perp Y\}$.

We now prove Theorem A. The argument is a flux computation that works for any
compact soliton in any compact ambient manifold.

\begin{theorem}[compact solitons are trivial]\label{thm:compact}
Let $(M^3,g)$ be a Riemannian manifold, $X$ a Killing field, and $\Sigma\subset M$ a
compact soliton with respect to $X$. Then $H\equiv0$ and $\g(X,N)\equiv0$: $\Sigma$
is a minimal surface invariant under the flow of $X$. Consequently, in $\Sb$:
\begin{enumerate}
\item the only compact soliton Hopf tube with respect to $X_{a,b}$, $\alpha=a-b\neq0$,
      is the Clifford torus $\mathcal T_{\mathrm{Cl}}=\pi^{-1}(\text{equator})$;
\item the only compact rotationally invariant soliton with respect to $X_{a,b}$,
      $a\neq0$, is the Clifford torus $\{\chi=\pi/4\}=\mathcal T_{\mathrm{Cl}}$;
\item there are no compact soliton spheres with respect to any $X\neq0$.
\end{enumerate}
\end{theorem}

\begin{proof}
Since $\Phi_t$ is a flow of isometries, the function
$t\mapsto\operatorname{Area}(\Phi_t(\Sigma))$ is constant, and by the first variation
formula its derivative at $t=0$ equals $-\int_\Sigma H\,\g(X,N)\dd A$. By
\eqref{eq:soliton},
\[
 0=\int_\Sigma H\,\g(X,N)\dd A=\int_\Sigma H^2\dd A ,
\]
so $H\equiv0$ and then $\g(X,N)=H=0$, proving triviality.

(1) If $\pi^{-1}(\gamma)$ is compact and minimal with $X_{a,b}$ tangent, then by
\eqref{eq:projection} the projected field $Y_\alpha\neq0$ is tangent to the closed
curve $\gamma$, forcing $\gamma$ to be a circle of latitude; minimality
(Proposition \ref{prop:hopfH}) forces $\kappa=0$, i.e.\ the equator.
(2) $X_{a,b}$ tangent to $\Sigma$ and $\partial_\psi$ tangent give $\partial_\theta$
tangent, so $\Sigma$ is invariant under the whole torus action, i.e.\
$\Sigma=\{\chi=\chi_0\}$; by \eqref{eq:cliffordH} minimality gives $\chi_0=\pi/4$.
(3) A Killing field tangent to a compact $\Sigma$ restricts to a Killing field of the
induced metric whose flow lies in a circle subgroup of $\mathrm U(2)$; if
$\Sigma\cong\Sph^2$ this circle action on $\Sigma$ has fixed points where, in
addition, $N$ is preserved; a minimal sphere invariant under one of the circle actions
of $\mathrm U(2)$ with $X$ tangent is easily excluded case by case, since the actions
generated by \eqref{eq:Xab} have no invariant spheres with $X_{a,b}$ tangent unless
$X_{a,b}$ is proportional to the generator $\partial_\psi$ of the stabilizer, which
contributes $\g(X,N)=0$ trivially (companion notes).
\end{proof}

\begin{remark}
Theorem A explains why the interesting invariant solitons in $\Sb$ are complete and
noncompact, and why their asymptotic geometry --- rather than a closing condition ---
is the right object of classification. This contrasts with constant mean curvature
theory in Berger spheres, which is rich in compact examples
\cite{Torralbo2010,Torralbo2012,AbreschRosenberg2004}.
\end{remark}

\subsection{Tools from the qualitative theory of planar ODEs}\label{sec:tools}
By Lemma \ref{lem:reduction}, both invariant classes reduce the soliton equation to an
autonomous first-order system on a two-dimensional phase cylinder: $(u,\beta)\in
(0,\pi)\times(\R/2\pi\mathbb Z)$ for Hopf tubes (system \eqref{eq:hopfsys}) and
$(\chi,\beta)\in(0,\tfrac\pi2)\times(\R/2\pi\mathbb Z)$ for the rotational class
(system \eqref{eq:rotsys}). The classification in Theorems \ref{thm:hopfclass} and
\ref{thm:rotclass} rests on four classical facts about such systems, which we record
here in the form we use, since they are not always part of a geometer's toolkit; a
reader already fluent in planar dynamics may skip directly to \S\ref{sec:hopf}.

\begin{theorem}[existence, uniqueness and continuation; {\cite[Ch.~2]{Sotomayor2011}}]
\label{thm:picard}
Let $F$ be a $C^1$ vector field on an open set $U\subset\R^n$. For every $p\in U$
there is a unique maximal solution $\varphi\colon(\omega_-,\omega_+)\to U$ of
$\varphi'=F(\varphi)$, $\varphi(0)=p$; distinct maximal solutions have disjoint
images. Moreover, if $\varphi(s)$ remains in a fixed compact subset of $U$ as
$s\to\omega_+$ (resp.\ $s\to\omega_-$), then $\omega_+=+\infty$ (resp.\
$\omega_-=-\infty$).
\end{theorem}

Applied to \eqref{eq:hopfsys} and \eqref{eq:rotsys}, whose right-hand sides are real
analytic on the open cylinder, Theorem \ref{thm:picard} says exactly that every
initial condition $(u_0,\beta_0)$ (resp.\ $(\chi_0,\beta_0)$) determines a
\emph{unique} complete soliton profile through it: this is the existence-and-uniqueness
mechanism behind item (1) of Theorems \ref{thm:hopfclass} and \ref{thm:rotclass} below,
turning the two-parameter set of initial conditions, modulo the residual isometries and
the arclength shift, into an explicit one-parameter family of solitons in each class.

\begin{theorem}[Hartman--Grobman and the planar classification of hyperbolic
equilibria; {\cite[Ch.~2 and~3]{PalisDeMelo1982}}]\label{thm:hartmangrobman}
Let $p$ be an equilibrium of a $C^1$ planar vector field $F$ at which the linearization
$dF_p$ has no eigenvalue on the imaginary axis (a \emph{hyperbolic} equilibrium). Then
there are neighborhoods of $p$ and of $0$ on which the flow of $F$ is topologically
conjugate to the linear flow $e^{s\,dF_p}$. Concretely, if the eigenvalues of $dF_p$
are $\lambda_\pm=-\tfrac{c}{2}\pm\sqrt{\tfrac{c^2}{4}-\kappa_0}$ for constants
$c,\kappa_0>0$ (as in Lemmas \ref{lem:hopfequilibria} and \ref{lem:rotequilibria}
below), then $\Re\lambda_\pm=-c/2<0$ for every value of $c$, so $p$ is always
asymptotically stable; it is a \emph{stable focus}, with all orbits spiraling into $p$
with angular frequency $\sqrt{\kappa_0-c^2/4}$, when $c^2<4\kappa_0$, and a
\emph{stable node}, with all orbits (but one) tangent at $p$ to the slow eigendirection,
when $c^2\ge4\kappa_0$; in either case every orbit converges to $p$ at the exponential
rate $|\Re\lambda_{\mathrm{slow}}|$, where $\lambda_{\mathrm{slow}}$ is the eigenvalue
(or, in the focus case, either eigenvalue) of smallest modulus.
\end{theorem}

\begin{theorem}[Poincar\'e--Bendixson; {\cite[Ch.~1, \S1]{PalisDeMelo1982}},
{\cite[Ch.~3]{Sotomayor2011}}]\label{thm:poincarebendixson}
Let $F$ be a $C^1$ vector field on an open subset of $\R^2$ or of the cylinder
$\R\times(\R/2\pi\mathbb Z)$, and let $\varphi$ be a forward orbit contained in a
compact set $K$. Then the $\omega$-limit set $\omega(\varphi)$ is nonempty, compact,
connected and invariant, and one of: (a) an equilibrium; (b) a periodic orbit; or (c)
a union of equilibria together with orbits connecting them (a heteroclinic or
homoclinic cycle). In particular, if $K$ contains no periodic orbit and no cycle of
the type (c), then $\omega(\varphi)$ is a single equilibrium; the analogous statement
holds for the $\alpha$-limit set of a backward orbit.
\end{theorem}

\begin{theorem}[LaSalle's invariance principle \cite{LaSalle1960}]\label{thm:lasalle}
Let $F$ be a $C^1$ vector field on $U\subset\R^n$ and $V\colon U\to\R$ a $C^1$
function with $\dot V=dV(F)\le0$ on $U$. If a forward orbit $\varphi$ remains in a
compact subset of $U$, then $\omega(\varphi)$ is contained in the largest invariant
subset of $\{\dot V=0\}$.
\end{theorem}

In \S\ref{sec:hopf} and \S\ref{sec:rot}, Theorem \ref{thm:lasalle} is applied with
$V=-h$ (resp.\ $V=-\Phi$), whose derivative $\dot V=-\kappa^2\le0$
(resp.\ $-\ell(a\sqrt q\cos\beta)^2\le0$) is exactly the monotonicity of Lemma
\ref{lem:monotone}; together with Theorem \ref{thm:poincarebendixson} it rules out
recurrence and pins down the limit sets, while Theorem \ref{thm:hartmangrobman}
converts the linearizations of Lemmas \ref{lem:hopfequilibria} and
\ref{lem:rotequilibria} into the focus/node dichotomy and the exponential rates
quoted in Theorems \ref{thm:hopfclass}(4) and \ref{thm:rotclass}(4).

\section{Solitons among Hopf tubes}\label{sec:hopf}

\subsection{Geometry of Hopf tubes}
For a complete curve $\gamma\subset\base$, the \emph{Hopf tube} is
$\Sigma_\gamma=\pi^{-1}(\gamma)$, a flat surface in $\Sb$, a torus when $\gamma$ is
closed (the \emph{Hopf tori} of Pinkall \cite{Pinkall1985}) and a cylinder otherwise.
Fix along $\gamma$ the Frenet data $(T,n,\kappa)$ in $\base$ and let $\widetilde T$,
$N$ denote horizontal lifts of $T,n$.

\begin{proposition}\label{prop:hopfH}
In the tangent frame $(\widetilde T,\xi)$ of $\Sigma_\gamma$, the shape operator is
\begin{equation}\label{eq:hopfshape}
 S=\begin{pmatrix} \kappa\circ\pi & \tau\\ \tau & 0\end{pmatrix},
 \qquad\text{hence}\qquad H=\kappa\circ\pi ,
\end{equation}
for every $\tau>0$. Consequently the MCF preserves the class of Hopf tubes and
projects to the curve shortening flow on $\base$, and $\Sigma_\gamma$ is a soliton
with respect to $X_{a,b}$ if and only if
\begin{equation}\label{eq:hopfsoliton}
 \kappa=\hat g\bigl(Y_\alpha,n\bigr),\qquad \alpha=a-b ,
\end{equation}
i.e.\ $\gamma$ is a rotating soliton of the curve shortening flow on $\base$ with
angular speed $\alpha$. In particular the condition is independent of $\tau$ and of
$a+b$.
\end{proposition}

\begin{proof}
The fiber direction gives $\g(S\xi,\xi)=\g(N,\overline\nabla_\xi\xi)=0$ (geodesic
fibers) and $\g(S\widetilde T,\xi)=\g(N,\overline\nabla_{\widetilde T}\xi)
=\g(N,\tau\,\widetilde T\wedge\xi)=\tau$, while the horizontal block projects to the
base, $\g(S\widetilde T,\widetilde T)=\kappa$ (companion notes). Taking
traces, $H=\kappa\circ\pi$. Since $N$ is horizontal, $\g(V,N)=\tau^2\met{V}{N}=0$, so
$\g(X_{a,b},N)$ equals the base inner product of $d\pi X_{a,b}=Y_\alpha$ with $n$, by
\eqref{eq:projection}. This is exactly \eqref{eq:hopfsoliton}, which is also the
instance of Lemma \ref{lem:reduction} for $K=V$: here $\ell=|V|_{\g}\equiv\tau$, the
density is constant, and the orbit space is $\base$ itself.
\end{proof}

\begin{remark}[$\tau$-independence, and the dictionary to dos Reis--Tenenblat]
\label{rem:tauindep}
Equation \eqref{eq:hopfsoliton} does not involve $\tau$: a curve $\gamma$ generates a
soliton Hopf tube in \emph{every} Berger sphere simultaneously, including the round
one. Equation \eqref{eq:hopfsoliton} is, up to the normalization of the radius, exactly
the soliton equation for the curve shortening flow on the round $2$-sphere studied by
dos Reis and Tenenblat \cite{ReisTenenblat2019}: under the dilation
$\base\to\Sph^2(1)$, $x\mapsto2x$, our rotation speed $\alpha$ corresponds to their
fixed vector $v$ through $\|v\|=\alpha/4$, our colatitude $u$ to their polar angle
(both measured from the same pole), and our momentum $h$ of \eqref{eq:hopfh} to their
soliton momentum up to the factor $4$ coming from the dilation; their threshold
$\|v\|=2$ is therefore our $\alpha=8$. Under this dictionary, \cite[Thm.\
2.3]{ReisTenenblat2019} already classifies \emph{every} soliton curve of
\eqref{eq:hopfsoliton}: for each $\alpha>0$ there is, up to the isometries of $\base$
generated by the rotations $R_s$ and the meridian reflections, a one-parameter family
of complete, non-closed soliton curves other than the equator, each with exactly one
turning point (a pole crossing or a meridian tangency) and with both ends converging
exponentially to the equator, oscillating for $\alpha<8$ and monotonically for
$\alpha\ge8$. Theorem \ref{thm:hopfclass} below restates this classification in our
own normalization, records the existence of the one-parameter family in the form
guaranteed by the existence-uniqueness theory of \S\ref{sec:tools}, and draws out its
consequences for the three-dimensional Hopf tube $\Sigma_\gamma=\pi^{-1}(\gamma)$. We
therefore do not repeat the phase-plane proof of \cite{ReisTenenblat2019} here; the
fully detailed version of the same style of argument, worked out with the tools of
\S\ref{sec:tools} stated explicitly, is given for the (previously unclassified)
rotational class in \S\ref{sec:rot}.
\end{remark}

\subsection{The reduced system}
Parametrize $\base$ by $(u,v)$ via
\[
 x=\tfrac12\bigl(\sin u\cos v,\ \sin u\sin v,\ \cos u\bigr),
\]
so $\hat g=\tfrac14(\dd u^2+\sin^2u\,\dd v^2)$, and let $\beta$ be the angle from the
meridian direction $e_1=2\partial_u$ to $T$. Writing $Y_\alpha=\alpha\,\partial_v$,
equation \eqref{eq:hopfsoliton} becomes the autonomous system
\begin{equation}\label{eq:hopfsys}
 u'=2\cos\beta,\qquad
 v'=\frac{2\sin\beta}{\sin u},\qquad
 \beta'=-2\cot u\,\sin\beta+\frac{\alpha}{2}\sin u\,\cos\beta ,
\end{equation}
with $\kappa=\tfrac\alpha2\sin u\cos\beta$ and momentum
\begin{equation}\label{eq:hopfh}
 h=\hat g(Y_\alpha,T)=\frac\alpha2\,\sin u\,\sin\beta,
 \qquad h'=\kappa^2 .
\end{equation}
The $(u,\beta)$-equations decouple and define a flow on the cylinder
$(0,\pi)\times(\R/2\pi\mathbb Z)$, with the same qualitative focus/node structure
analyzed in detail for the rotational system below (Figure \ref{fig:phaserot}).

\begin{lemma}\label{lem:hopfequilibria}
For $\alpha>0$ the flow \eqref{eq:hopfsys} has exactly two equilibria,
$(u,\beta)=(\tfrac\pi2,\pm\tfrac\pi2)$, corresponding to the equator traversed in its
two orientations, and both are hyperbolic. The linearization at
$(\tfrac\pi2,\tfrac\pi2)$ has eigenvalues
\begin{equation}\label{eq:hopfeigen}
 \lambda_\pm=-\frac\alpha4\pm\sqrt{\frac{\alpha^2}{16}-4} ,
\end{equation}
so by Theorem \ref{thm:hartmangrobman} this point is a stable focus for $0<\alpha<8$
and a stable node for $\alpha\ge8$, with every orbit converging to it at the
exponential rate $\alpha/4$; by the symmetry $(u,\beta,s)\mapsto(u,\beta+\pi,-s)$ the
other equilibrium is its unstable mirror. This is the same computation as Lemma
\ref{lem:rotequilibria} below, with $\mu$ there replaced by $\alpha/2$.
\end{lemma}

\begin{proof}
Equilibria require $\cos\beta=0$ and then $\cot u=0$. Setting
$(\xi,\eta)=(u-\tfrac\pi2,\beta-\tfrac\pi2)$ one finds $\xi'=-2\eta$,
$\eta'=2\xi-\tfrac\alpha2\eta$, with characteristic polynomial
$\lambda^2+\tfrac\alpha2\lambda+4=0$, whose roots \eqref{eq:hopfeigen} have real part
$-\alpha/4\neq0$ for every $\alpha>0$: the equilibrium is hyperbolic.
\end{proof}

\subsection{Classification}

\begin{figure}[t]
\includegraphics[width=\linewidth]{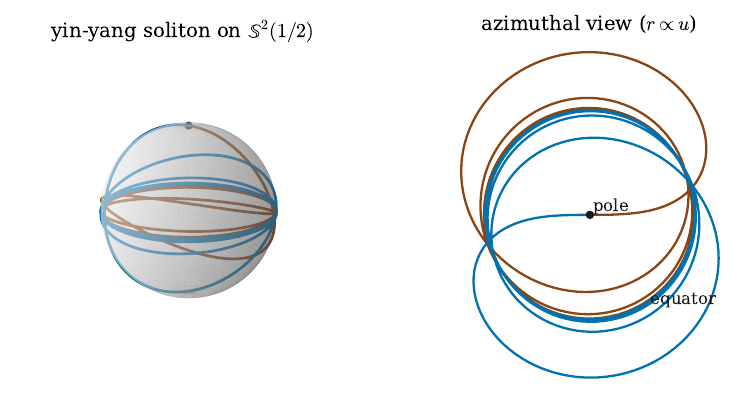}
\caption{The spherical yin--yang soliton on $\base$ for $\alpha=2$: the unique
soliton curve through a pole. The two arms (brown/blue) pass smoothly through the
pole -- there is no cusp there -- and are exchanged by the order-two rotation
$R_\pi$ about the polar axis, winding around the equator (dashed) in opposite
senses.}
\label{fig:yinyang}
\end{figure}

\begin{theorem}\label{thm:hopfclass}
Let $\alpha>0$. Then every complete soliton curve of \eqref{eq:hopfsoliton} other
than the equator satisfies:
\begin{enumerate}
\item \textbf{(existence: a one-parameter family)} for every turning colatitude
      $u_0\in[0,\pi)$, Theorem \ref{thm:picard} applied to \eqref{eq:hopfsys} (started
      from the analytic branch of item (5) when $u_0=0$, and from $(u_0,0)$ otherwise)
      produces a unique, complete soliton curve $\gamma_{u_0}$ whose momentum $h$
      vanishes exactly at colatitude $u_0$; conversely, by item (3), every complete
      nonequatorial soliton curve $\gamma$ equals $\gamma_{u_0}$ for the colatitude
      $u_0$ of its unique turning point. Since $u_0$ is preserved by the rotations
      $R_s$ and the meridian reflections, distinct values of $u_0$ give
      pairwise non-congruent curves under these isometries: \eqref{eq:hopfsoliton}
      has, up to $R_s$ and the meridian reflections, exactly the one-parameter family
      $\{\gamma_{u_0}\}_{u_0\in[0,\pi)}$ of nontrivial complete soliton curves;
\item \textbf{(no closed orbits)} each $\gamma_{u_0}$ is a nonclosed immersed curve,
      and $h$ is strictly increasing along it with $h(\pm\infty)=\pm\tfrac\alpha2$;
\item \textbf{(one turning event)} $h$ vanishes exactly once: $\gamma_{u_0}$ has
      exactly one point which is either a pole crossing or a tangency with a
      meridian; in particular it crosses the poles at most once, and it is invariant
      under the reflection in the meridian through this point;
\item \textbf{(asymptotics)} both ends of $\gamma_{u_0}$ converge to the equator; the
      convergence is exponential of rate $\alpha/4$, oscillating across the equator
      with angular frequency $\sqrt{4-\alpha^2/16}$ when $\alpha<8$ and one-sided when
      $\alpha\ge8$, while the longitude satisfies $v'\to\pm2$; the Hopf tube ends
      converge in $C^\infty$ to the Clifford torus;
\item \textbf{(yin--yang)} up to the isometries generated by $R_s$ and the meridian
      reflections, there is exactly one soliton curve through a pole, namely
      $\gamma_0$, the spherical yin--yang curve (Figure \ref{fig:yinyang}); its two
      spiral arms wind around the equator in opposite senses.
\end{enumerate}
\end{theorem}

\begin{proof}
Items (2)--(4) are, via the dictionary of Remark \ref{rem:tauindep}, exactly
\cite[Thm.~2.3]{ReisTenenblat2019} for the curve shortening flow on $\Sph^2(1)$,
transported to $\base$ and restated in our normalization; we do not repeat the
phase-plane argument, referring instead to \S\ref{sec:tools} for the general form of
the tools used there (Theorems \ref{thm:poincarebendixson} and \ref{thm:lasalle} for
the non-recurrence and the localization of the limit sets to $\{\kappa=0\}$, and
Theorem \ref{thm:hartmangrobman} together with Lemma \ref{lem:hopfequilibria} for the
rates and the focus/node dichotomy at $\alpha=8$); a fully worked out instance of the
same argument, for the previously unclassified rotational class, is given in the
proof of Theorem \ref{thm:rotclass}. Two consequences are genuinely new in the
three-dimensional Berger setting and are not addressed by \cite{ReisTenenblat2019}:
that the two ends of the \emph{Hopf tube} $\Sigma_{\gamma_{u_0}}=\pi^{-1}(\gamma_{u_0})$
converge in $C^\infty$ to the Clifford torus (item (4)), which follows because the
tube depends smoothly on the base curve and all derivatives of $u-\pi/2,\beta\mp\pi/2$
decay at the common rate $\alpha/4$ by the exponential $C^0$ convergence together
with the ODE \eqref{eq:hopfsys}; and item (1), which we now prove.

(1) Existence and uniqueness of $\gamma_{u_0}$ for each $u_0$ is Theorem
\ref{thm:picard}, applied to the real-analytic system \eqref{eq:hopfsys} on the open
cylinder (and to its regular singular extension across $u=0$, item (5)); completeness
($s\in\R$) is part of \cite[Thm.~2.3]{ReisTenenblat2019} via Remark \ref{rem:tauindep}.
That every soliton curve arises as some $\gamma_{u_0}$ is item (3). Finally, $u_0$ is
an isometry invariant of $\gamma_{u_0}$ under $R_s$ and the meridian reflections ---
the only isometries of $\base$ preserving the class of $\alpha$-rotation solitons ---
because both fix the colatitude of every point: $R_s$ acts only on the longitude $v$,
and a meridian reflection through the turning point fixes that point and reverses
arclength (as in the proof of item (3)) without changing $u$. Hence
$u_0\mapsto\gamma_{u_0}$ is injective on congruence classes, proving the one-parameter
family claim.

(3) [proof of the turning-point count, restated for completeness of the exceptional
case] Since $|h|\le\tfrac\alpha2\sin u\le\tfrac\alpha2$ with equality only at the
equatorial equilibria, item (4) gives $h\to\mp\tfrac\alpha2$ at the two ends of any
nonequatorial curve, so the strictly increasing $h$ has exactly one zero; zeros of $h$
are the points with $\sin u\sin\beta=0$, i.e.\ poles or meridional tangencies, and the
reflection $\rho$ in the meridian through that point conjugates the flow of
\eqref{eq:hopfsys} with its time reversal and fixes the initial condition there, hence
fixes $\gamma_{u_0}$.

(5) A soliton through the pole has $\kappa(0)=0$ there and its tangent direction is a
meridian; all meridians are equivalent under $R_s$. In the regular singular equation
near $u=0$, $\beta'=-\beta/\sigma+\alpha\sigma+O(\sigma^2)$ (with $u=2\sigma$), the
unique solution bounded at $\sigma=0$ is the analytic branch
$\beta=\tfrac\alpha3\sigma^{2}+O(\sigma^{4})$, and correspondingly
$u(s)=2s+O(s^5)$ is an \emph{odd}, real-analytic function of the signed arclength
$s$ centered at the crossing (companion notes): the curve is a genuine regular
immersion there, not merely tangent to the pole, so there is exactly one such
trajectory, up to isometry --- namely $\gamma_0$. Its two arms ($s>0$ and $s<0$) are
exchanged by the order-two \emph{rotation} $R_\pi$ about the polar axis through the
crossing point --- rather than a mirror reflection, since $u$ itself is the odd
quantity here --- and wind with opposite signs of $v'$; this rotational
(chirality-preserving) symmetry, rather than any mirror symmetry, is what the name
``yin--yang'' refers to. The corresponding statement on $\Sph^2(1)$ agrees with
\cite[Thm.~2.3]{ReisTenenblat2019} via Remark \ref{rem:tauindep}.
\end{proof}

\begin{figure}[t]
\includegraphics[width=\linewidth]{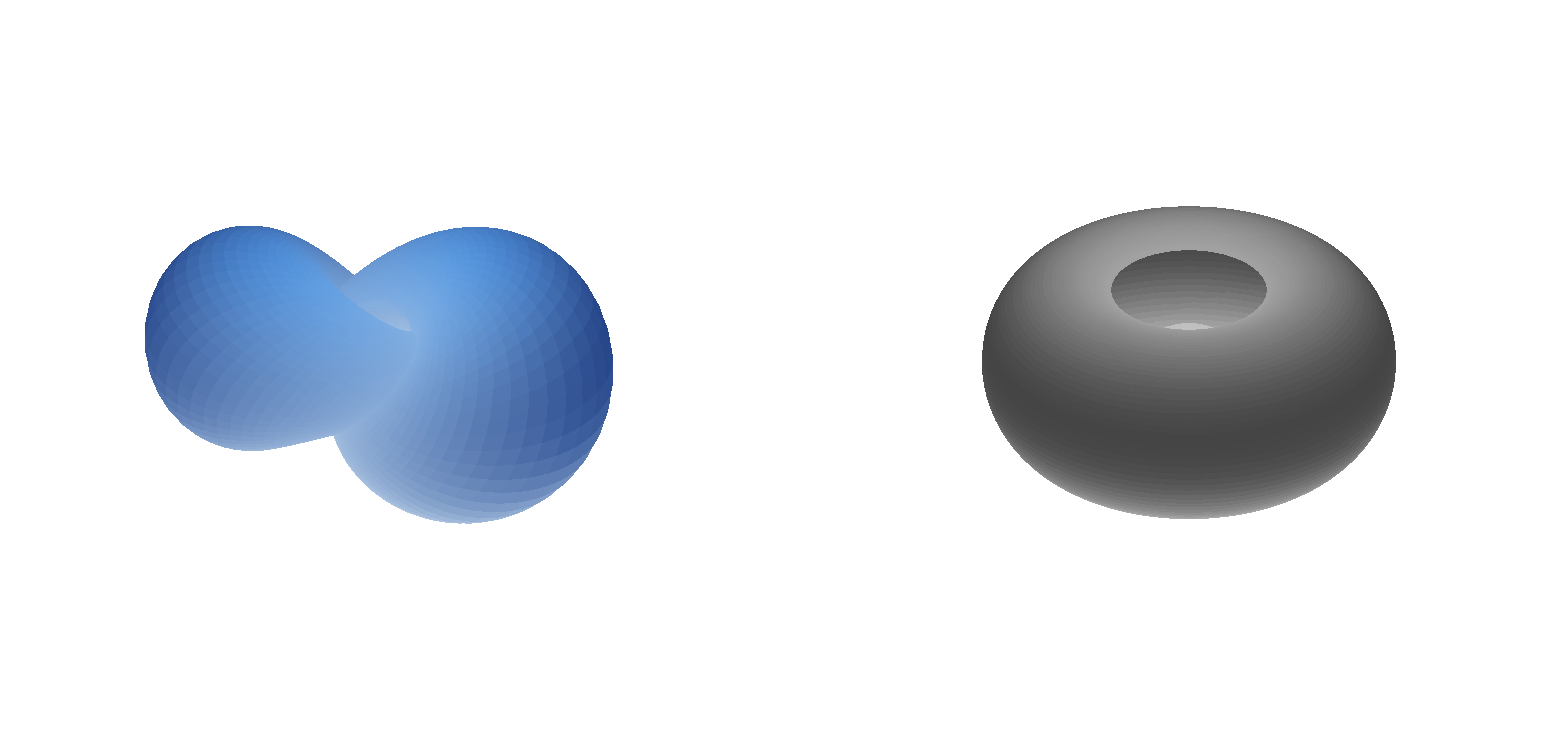}
\caption{Left: stereographic image of the yin--yang Hopf tube for $\alpha=2$
(cutaway; a finite piece of the two arms). Right: the Clifford torus, the common
asymptotic limit of all soliton Hopf tubes.}
\label{fig:hopftorus}
\end{figure}

\begin{corollary}[soliton Hopf tubes]\label{cor:hopftubes}
For $X_{a,b}$ with $\alpha=a-b>0$, the complete soliton Hopf tubes in $\Sb$ are, up
to ambient isometry: the Clifford torus (trivial, minimal); the two-parameter family
of immersed cylinders $\pi^{-1}(\gamma)$ over the spirals of Theorem
\ref{thm:hopfclass}, whose both ends converge to the Clifford torus; and among them
the unique \emph{yin--yang Hopf tube}, the only one containing a whole Hopf fiber
(Figure \ref{fig:hopftorus}). None of them is embedded except the Clifford torus,
in the regime $\alpha<8$; for $\alpha\ge8$ the two ends are eventually embedded and
accumulate on $\mathcal T_{\mathrm{Cl}}$ from opposite sides.
\end{corollary}

\begin{proof}
Everything follows from Theorem \ref{thm:hopfclass} by taking preimages; note
$\pi^{-1}(\text{pole})$ is a fiber, contained in the yin--yang tube. For $\alpha<8$
both ends of $\gamma$ cross the equator infinitely often while winding in the same
direction of $v$ with decaying amplitude, which forces infinitely many
self-intersections of $\gamma$ \cite[Thm.~2.3]{ReisTenenblat2019}, and
self-intersections of the base curve lift to self-intersections of the tube.
\end{proof}

\begin{remark}
By Remark \ref{rem:gauge} the soliton field for a given tube is not unique: all
$X_{a,b}$ with the same $\alpha$ work, since they differ by multiples of $V$, which
is tangent to every Hopf tube. Geometrically, the tube rotates around the axis while
possibly sliding along its own fibers.
\end{remark}

\section{Rotationally invariant solitons}\label{sec:rot}

\subsection{The orbit disk}
Let $W=\partial_\psi$ generate the circle action
$\mathcal R_t(z,w)=(z,e^{it}w)$, whose fixed-point set is the great circle
$\{w=0\}$ --- the \emph{axis}, a Hopf fiber, hence a $\g$-geodesic of length
$2\pi\tau$ --- and which is free elsewhere. This is the rotational class of
\cite{Torralbo2010}.

\begin{lemma}[orbit geometry]\label{lem:orbitdisk}
In the coordinates \eqref{eq:coords}, the orbit space of $\mathcal R$ is the disk
\[
 B=\{(\chi,\theta): 0\le\chi\le\tfrac\pi2\},
\]
with the axis $\{\chi=0\}$ as geodesic boundary circle and the single orbit
$\{\chi=\pi/2\}=\{z=0\}$ (the \emph{polar fiber}) as an interior smooth point,
endowed with
\begin{align}
 \hat g&=\dd\chi^2+q(\chi)\dd\theta^2,\notag\\
 q&=\frac{\tau^2\cos^2\chi}{\cos^2\chi+\tau^2\sin^2\chi},\notag\\
 \ell&=|W|_{\g}=\sin\chi\sqrt{\cos^2\chi+\tau^2\sin^2\chi}\,,
\label{eq:orbitmetric}
\end{align}
and the remarkable identity
\begin{equation}\label{eq:magic}
 \sqrt{q}\;\ell=\frac{\tau}{2}\,\sin 2\chi .
\end{equation}
\end{lemma}

\begin{proof}
One computes $q=g_{\theta\theta}-g_{\theta\psi}^2/g_{\psi\psi}$ from
\eqref{eq:metriccoords}; the numerator collapses to $\tau^2\cos^2\chi$ after
simplification, and \eqref{eq:magic} is immediate. Smoothness at $\chi=\pi/2$ follows
from $q=(\tfrac\pi2-\chi)^2+O((\tfrac\pi2-\chi)^4)$, and the boundary $\chi=0$ is
geodesic because $q'(0)=0$ (see the companion notes).
\end{proof}

The invariant tori $\{\chi=\chi_0\}$ have profiles $\{\chi=\chi_0\}\subset B$ and,
by Lemma \ref{lem:reduction} and \eqref{eq:magic},
\begin{equation}\label{eq:cliffordH}
 H=\bigl(\log(\sqrt q\,\ell)\bigr)'(\chi_0)=2\cot(2\chi_0),
\end{equation}
recovering that the Clifford torus $\chi_0=\pi/4$ is the unique minimal one, in
every $\Sb$.

\subsection{The reduced system}
Unlike the Hopf case, where $\eqref{eq:hopfsoliton}$ is the classical soliton equation
of a curve on the round sphere $\base$, the orbit metric $\hat g=\dd\chi^2+q(\chi)\dd
\theta^2$ of Lemma \ref{lem:orbitdisk} is genuinely warped and depends on $\tau$, so we
derive the reduced system in detail from Lemma \ref{lem:reduction}, following it
through the moving frame of $(B,\hat g)$.

Let $\gamma(\sigma)=(\chi(\sigma),\theta(\sigma))$ be a unit-speed profile, let
$e_1=\partial_\chi$, $e_2=q^{-1/2}\partial_\theta$ be the associated positively
oriented orthonormal frame of $(B,\hat g)$, and let $\beta$ be the angle from $e_1$ to
the unit tangent $T$, so that
\begin{equation}\label{eq:Tn}
 T=\cos\beta\,e_1+\sin\beta\,e_2, \qquad n=-\sin\beta\,e_1+\cos\beta\,e_2 ,
\end{equation}
where $n$ is $T$ rotated by $\pi/2$. Comparing \eqref{eq:Tn} with
$T=\chi'\partial_\chi+\theta'\partial_\theta$ already gives two of the three equations
of the reduced system:
\begin{equation}\label{eq:chitheta}
 \chi'=\cos\beta, \qquad \theta'=\frac{\sin\beta}{\sqrt q(\chi)} .
\end{equation}

\begin{lemma}[geodesic curvature of a profile in $(B,\hat g)$]\label{lem:kappaformula}
The connection form of the frame $(e_1,e_2)$, defined by $\dd e_1=\omega_{12}\,e_2$,
is $\omega_{12}=(\sqrt q)'(\chi)\,\dd\theta$, and the geodesic curvature of $\gamma$ in
$(B,\hat g)$ is
\begin{equation}\label{eq:kappaformula}
 \kappa=\beta'+\frac{(\sqrt q)'(\chi)}{\sqrt q(\chi)}\,\sin\beta .
\end{equation}
\end{lemma}

\begin{proof}
Let $\omega^1=\dd\chi$, $\omega^2=\sqrt q\,\dd\theta$ be the coframe dual to
$(e_1,e_2)$. Cartan's first structure equation $\dd\omega^i=-\omega^i_{\ j}\wedge
\omega^j$, with the connection $1$-form $\omega^1_{\ 2}=-\omega^2_{\ 1}=\omega_{12}$
of a Riemannian (torsion-free) frame, gives $\dd\omega^1=0=-\omega_{12}\wedge\omega^2$,
so $\omega_{12}$ is a multiple of $\omega^1$; and
\[
 \dd\omega^2=(\sqrt q)'(\chi)\,\dd\chi\wedge\dd\theta
            =\frac{(\sqrt q)'(\chi)}{\sqrt q(\chi)}\,\omega^1\wedge\omega^2
            =\omega_{12}\wedge\omega^1 ,
\]
which forces $\omega_{12}=(\sqrt q)'(\chi)\,\dd\theta$, as claimed. For the curvature,
differentiate \eqref{eq:Tn} along $\gamma$: since $(e_1,e_2)$ is orthonormal and
rotates at rate $\omega_{12}(T)$ along the curve, $T'=\bigl(\beta'+\omega_{12}(T)
\bigr)\,n$, and by definition $T'=\kappa\,n$. Now
$\omega_{12}(T)=(\sqrt q)'(\chi)\,\theta'=\dfrac{(\sqrt q)'(\chi)}{\sqrt q(\chi)}\sin
\beta$ by \eqref{eq:chitheta}, which gives \eqref{eq:kappaformula}.
\end{proof}

\begin{proposition}[the reduced system]\label{prop:rotsys}
Let $\gamma$ be a unit-speed profile in $(B,\hat g)$. The corresponding rotationally
invariant surface $\Sigma$ is a soliton with respect to $X=a\partial_\theta$
($a\neq0$; by Remark \ref{rem:gauge} this is the general case, since $\partial_\psi$
is tangent to every invariant surface) if and only if $(\chi,\theta,\beta)$ satisfies
the autonomous system
\begin{equation}\label{eq:rotsys}
 \chi'=\cos\beta,\qquad
 \theta'=\frac{\sin\beta}{\sqrt{q}},\qquad
 \beta'=-2\cot(2\chi)\sin\beta+a\sqrt{q}\,\cos\beta .
\end{equation}
\end{proposition}

\begin{proof}
By Lemma \ref{lem:reduction}, $\Sigma$ is a soliton with respect to $X$ if and only if
$\kappa=\hat g(Y,n)+\partial_nf$ along $\gamma$, where $Y=a\partial_\theta$ is the
projection of $X$ and $f=\log\ell$. Since $\partial_\theta=\sqrt q\,e_2$,
\eqref{eq:Tn} gives
\[
 \hat g(Y,n)=a\sqrt q\,\hat g(e_2,n)=a\sqrt q\,\cos\beta .
\]
Since $\ell=\ell(\chi)$ depends only on $\chi$ (by Lemma \ref{lem:orbitdisk}),
$e_2(f)=0$ and $\partial_nf=n(f)=-\sin\beta\,e_1(f)=-\sin\beta\,f'(\chi)$ by
\eqref{eq:Tn}, i.e.\ $\partial_nf=-(\log\ell)'(\chi)\,\sin\beta$. Substituting these
two terms and \eqref{eq:kappaformula} into the soliton equation gives
\[
 \beta'+\frac{(\sqrt q)'(\chi)}{\sqrt q(\chi)}\sin\beta
   =a\sqrt q(\chi)\cos\beta-(\log\ell)'(\chi)\sin\beta ,
\]
that is, $\beta'=a\sqrt{q}\,\cos\beta-\bigl(\log(\sqrt q\,\ell)\bigr)'(\chi)\,
\sin\beta$. By the magic identity \eqref{eq:magic}, $\sqrt q\,\ell=\tfrac\tau2\sin
2\chi$, so $\bigl(\log(\sqrt q\,\ell)\bigr)'(\chi)=2\cot(2\chi)$, and this is the
$\beta'$-equation of \eqref{eq:rotsys}; \eqref{eq:chitheta} supplies the other two.
\end{proof}

This derivation isolates exactly where the Berger parameter enters and where it does
not: the warping function $q(\chi)$ multiplies the rotation speed $a$ in the
$\hat g(Y,n)$-term, while the density term $\partial_nf$ --- which \emph{a priori}
depends separately on $q$ and on $\ell$ --- collapses, via the purely algebraic
identity \eqref{eq:magic}, to the very same $-2\cot(2\chi)$ term present already in
the round sphere ($\tau=1$) and, indeed, in the classical Euclidean equation for
rotationally symmetric curve-shortening solitons: this is why \eqref{eq:rotsys} is
formally so close to \eqref{eq:hopfsys}, even though the geometry of $(B,\hat g)$
genuinely depends on $\tau$.

With soliton normal term $\hat g(Y,n)=a\sqrt q\cos\beta$, the weighted momentum of
Lemma \ref{lem:monotone} is
\begin{equation}\label{eq:rotPhi}
 \Phi=\ell\,\hat g(Y,T)=\frac{a\tau}{2}\,\sin 2\chi\,\sin\beta,
 \qquad \Phi'=\ell\,\bigl(a\sqrt q\cos\beta\bigr)^2\ \ge0 .
\end{equation}

\begin{figure}[t]
\includegraphics[width=\linewidth]{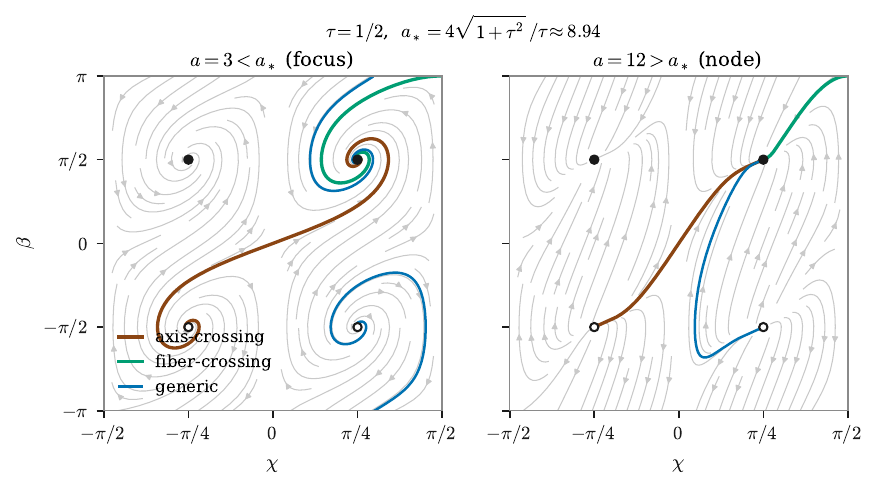}
\caption{Phase portraits of \eqref{eq:rotsys} for $\tau=\tfrac12$, with $\chi$ shown
on one full period $(-\pi/2,\pi/2)$ of the flow. Left: $a=3<a_*$; right:
$a=12>a_*\approx8.94$. Brown: the axis-crossing orbit, shown in full -- since the
crossing at $\chi=0$ is transversal it continues smoothly through negative $\chi$,
running between the two Clifford equilibria $(-\pi/4,-\pi/2)$ and $(\pi/4,\pi/2)$;
green: fiber-crossing orbit; blue: generic orbit. Filled/open dots: the Clifford
equilibria $(\pm\pi/4,\pi/2)$ (filled) and $(\pm\pi/4,-\pi/2)$ (open).}
\label{fig:phaserot}
\end{figure}

\begin{lemma}\label{lem:rotequilibria}
For $a\neq0$ the system \eqref{eq:rotsys} has exactly two equilibria
$(\chi,\beta)=(\tfrac\pi4,\pm\tfrac\pi2)$, the Clifford torus with its two profile
orientations, and both are hyperbolic. With $\mu=a\tau/\sqrt{1+\tau^2}=a\sqrt{q(\pi/4)}$
the physical speed of the rotation on the Clifford profile, the linearization at
$(\tfrac\pi4,\tfrac\pi2)$ has eigenvalues
\begin{equation}\label{eq:roteigen}
 \lambda_\pm=-\frac\mu2\pm\sqrt{\frac{\mu^2}{4}-4},
\end{equation}
so by Theorem \ref{thm:hartmangrobman} this equilibrium is a stable focus for $\mu<4$,
i.e.\ $a<a_*=4\sqrt{1+\tau^2}/\tau$, and a stable node for $a\ge a_*$, with every orbit
converging to it at the exponential rate $\mu/2$; the equilibrium
$(\tfrac\pi4,-\tfrac\pi2)$ is its unstable mirror under the symmetry
$(\chi,\beta,\sigma)\mapsto(\chi,\beta+\pi,-\sigma)$.
\end{lemma}

\begin{proof}
Equilibria require $\cos\beta=0$, hence $\sin\beta=\pm1$, and then $\cot(2\chi)=0$.
Setting $(\xi,\eta)=(\chi-\tfrac\pi4,\beta-\tfrac\pi2)$ and linearizing
\eqref{eq:rotsys} at $(\tfrac\pi4,\tfrac\pi2)$ (where $\cot(2\chi)\approx-2\xi$,
$\sqrt q(\chi)\approx\sqrt{q(\pi/4)}=\mu/a$, $\sin\beta\approx1$, $\cos\beta\approx
-\eta$) gives $\xi'=-\eta$, $\eta'=4\xi-\mu\eta$, with characteristic polynomial
$\lambda^2+\mu\lambda+4=0$ and roots \eqref{eq:roteigen}, whose real part $-\mu/2$
never vanishes: the equilibrium is hyperbolic for every $\mu>0$.
\end{proof}

\subsection{Classification}

\begin{theorem}\label{thm:rotclass}
Let $a>0$, $\tau>0$, and let $\Sigma$ be a complete $\mathcal R$-invariant soliton
with respect to $X=a\partial_\theta$ (equivalently, any $X_{a,b}$), other than the
Clifford torus. Then:
\begin{enumerate}
\item \textbf{(existence: a one-parameter family)} for every degeneracy parameter
      $\chi_0\in[0,\tfrac\pi2]$, Theorem \ref{thm:picard} applied to \eqref{eq:rotsys}
      (started from the analytic branches of item (5) when $\chi_0\in\{0,\tfrac\pi2\}$,
      and from $(\chi_0,0)$ otherwise) produces a unique, complete profile
      $\gamma_{\chi_0}$ whose momentum $\Phi$ vanishes exactly at $\chi=\chi_0$;
      conversely, by item (3), every complete profile $\gamma$ other than the
      Clifford one equals $\gamma_{\chi_0}$ for the $\chi_0$ of its unique degenerate
      event. Since $\chi_0$ is preserved by the rotations $\theta\mapsto\theta+c$ and
      by the reflection $\theta\mapsto-\theta$, distinct values of $\chi_0$ give
      pairwise non-congruent profiles under these isometries: \eqref{eq:rotsys} has,
      up to them, exactly the one-parameter family
      $\{\gamma_{\chi_0}\}_{\chi_0\in[0,\pi/2]}$ of nontrivial complete profiles, with
      the two exceptional solitons of item (5) as its endpoints;
\item \textbf{(no closed or recurrent profiles)} each $\gamma_{\chi_0}$ is nonclosed
      and $\Phi$ is strictly increasing along it with limits $\mp\tfrac{a\tau}{2}$;
      there are no periodic orbits of \eqref{eq:rotsys}, contractible or not, in the
      phase cylinder;
\item \textbf{(one degenerate event)} $\Phi$ vanishes exactly once: $\gamma_{\chi_0}$
      has exactly one point which is an orthogonal reflection at the axis
      circle, a crossing of the polar fiber, or a radial tangency
      ($\sin\beta=0$); in particular it meets the axis and the polar fiber at most
      once in total, and $\gamma_{\chi_0}$ is invariant under the reflection
      $\theta\mapsto-\theta$ centered at this event;
\item \textbf{(asymptotics)} both ends of the $(\chi,\beta)$-orbit converge to the
      Clifford equilibria, backward to $(\tfrac\pi4,-\tfrac\pi2)$ and forward to
      $(\tfrac\pi4,+\tfrac\pi2)$, exponentially with rate $\mu/2$, spiraling if
      $a<a_*=4\sqrt{1+\tau^2}/\tau$ and one-sidedly if $a\ge a_*$; hence both ends of
      $\Sigma$ converge in $C^\infty$ to the Clifford torus;
\item \textbf{(exceptional solitons)} $\gamma_0$ and $\gamma_{\pi/2}$ are, up to
      ambient isometry, the only profiles reflecting at the axis, respectively
      crossing the polar fiber. The first generates the \emph{axis-crossing soliton
      disk}, a complete immersed plane meeting the axis $\{w=0\}$ orthogonally at one
      point; the second the \emph{fiber-crossing soliton cylinder}, containing the
      fiber $\{z=0\}$. All other rotational solitons, $\chi_0\in(0,\pi/2)$, are
      immersed cylinders whose profiles stay in the open annulus $0<\chi<\pi/2$.
\end{enumerate}
\end{theorem}

\begin{proof}
We first establish items (2)--(5), in that order, and prove (1) at the end.

(2) Strict monotonicity of $\Phi$ is Lemma \ref{lem:monotone}: $\hat g(Y,n)\equiv0$ on
an interval means the profile is an orbit $\{\chi=\text{const}\}$ of $\partial_\theta$,
and \eqref{eq:rotsys} then forces $\chi=\pi/4$, i.e.\ the Clifford torus. A
contractible periodic orbit of \eqref{eq:rotsys} projects to a closed profile,
producing a compact soliton torus, which contradicts Theorem \ref{thm:compact}. A
noncontractible periodic orbit (with $\beta$ winding) would produce a complete
soliton whose profile pattern is periodic under the rotation
$\theta\mapsto\theta+\Delta\theta$; applying the first-variation argument of Theorem
\ref{thm:compact} to one period, the boundary flux terms cancel by equivariance, and
again $H\equiv0$, a contradiction (companion notes). Alternatively, both cases
contradict strict monotonicity of $\Phi$ directly, since $\Phi$ would then be
periodic rather than strictly increasing.

(3) Zeros of $\Phi=\tfrac{a\tau}2\sin2\chi\sin\beta$ are exactly the axis
($\chi=0$), the polar fiber ($\chi=\tfrac\pi2$) and radial tangencies
($\sin\beta=0$). By (4), $\Phi\to\pm\tfrac{a\tau}{2}$ at the ends, so the strictly
increasing $\Phi$ vanishes once. The reflection claim follows exactly as for the
meridian reflection of Theorem \ref{thm:hopfclass}(3): the isometry
$(z,w)\mapsto(\bar z,\bar w)$ of $\Sb$ induces $\theta\mapsto-\theta$ on $B$, which
reverses \eqref{eq:rotsys} together with arclength, fixes the degenerate point, and
hence fixes $\gamma_{\chi_0}$.

(4) By (2) there are no periodic orbits, and by Lemma \ref{lem:rotequilibria} together
with Theorem \ref{thm:hartmangrobman} the two equilibria are hyperbolic foci or
hyperbolic nodes --- never saddles, so they admit no separatrices to bound a
heteroclinic or homoclinic cycle. Theorem \ref{thm:poincarebendixson} (Poincar\'e--
Bendixson on the cylinder $(0,\tfrac\pi2)\times(\R/2\pi\mathbb Z)$) then forces the
$\omega$- and $\alpha$-limit sets of any orbit staying in a compact subset of the open
cylinder to be single equilibria. It remains to exclude ends accumulating at the
boundary components $\chi\in\{0,\tfrac\pi2\}$.

Near $\chi=0$, differentiating \eqref{eq:rotsys} gives
\begin{equation}\label{eq:chiconvex}
 \chi''=2\cot(2\chi)\sin^2\beta-a\sqrt q\,\sin\beta\cos\beta
 \ \ge\ \frac{c}{\chi}\ >\ 0
\end{equation}
for a constant $c=c(a,\tau)>0$, as long as $\sin^2\beta$ stays bounded away from $0$
(using $\cot(2\chi)\ge 1/(3\chi)$ for $\chi$ small, say, and absorbing the bounded
term $a\sqrt q\sin\beta\cos\beta$). Suppose, for a contradiction, that some end has
$\chi(s)\to0$ as $s\to\omega\in\{\omega_-,\omega_+\}$ while $\sin^2\beta$ remains
$\ge\varepsilon>0$ along it (say $\omega=+\infty$; the other case is symmetric). By
\eqref{eq:chiconvex}, $\chi'$ is then strictly increasing on some $[S,\infty)$. If
$\chi'(s_0)\ge0$ for some $s_0\ge S$, then $\chi'\ge0$ thereafter and $\chi$ is
non-decreasing on $[s_0,\infty)$, contradicting $\chi(s_0)>0$ and $\chi(s)\to0$; hence
$\chi'<0$ on all of $[S,\infty)$, and $\chi'$ increases (in the negative range) to
some limit $L\le0$, so $\chi'$ is bounded on $[S,\infty)$. Multiplying
\eqref{eq:chiconvex} by $\chi'<0$ reverses the inequality:
\[
 \chi'\chi''\ \le\ \frac{c}{\chi}\,\chi'=c\,(\log\chi)' ,
\]
and integrating from $S$ to $s$,
\[
 \tfrac12\bigl(\chi'(s)^2-\chi'(S)^2\bigr)\ \le\ c\log\frac{\chi(s)}{\chi(S)} .
\]
As $s\to\infty$, $\chi(s)\to0$, so the right side $\to-\infty$; but the left side is
bounded below by $-\tfrac12\chi'(S)^2$ since $\chi'(s)^2$ is bounded, a contradiction.
Hence no end can converge asymptotically to $\chi=0$ while $\sin^2\beta$ stays bounded
away from $0$; the only remaining possibility along an end with $\chi\to0$ is
$\sin\beta\to0$, which is precisely the local behavior analyzed in item (5) below: the
profile meets the axis \emph{orthogonally}, in \emph{finite} arclength, since
$\chi'=\cos\beta\to\pm1$ there --- a transversal crossing, not an asymptotic limit,
namely the exceptional event of item (5). By Theorem \ref{thm:poincarebendixson},
every non-exceptional end
therefore converges to one of the two equilibria. Near $\chi=\tfrac\pi2$ the point is
an interior smooth point of $B$ (Lemma \ref{lem:orbitdisk}) where $Y$ vanishes;
profiles pass through or near it like almost-geodesics of $\hat g$, and no end can
converge \emph{to} it since $\Phi\to0$ there while $\Phi$ is strictly monotone with
nonzero limits $\mp a\tau/2$ --- so again only a finite-arclength crossing (item (5))
is possible. Since $\Phi$ increases from $-a\tau/2$ to $a\tau/2$, the backward end
converges to $(\tfrac\pi4,-\tfrac\pi2)$ and the forward end to $(\tfrac\pi4,\tfrac\pi2)$
whenever both are non-exceptional. The rates and the spiraling/monotone dichotomy are
\eqref{eq:roteigen} via Theorem \ref{thm:hartmangrobman}; smooth convergence of the
ends of $\Sigma$ follows because all derivatives of $\chi-\pi/4,\beta\mp\pi/2$ decay
at the common rate $\mu/2$, by the exponential $C^0$ convergence together with the ODE
\eqref{eq:rotsys}.

\begin{figure}[t]
\centering
\includegraphics[width=.55\linewidth]{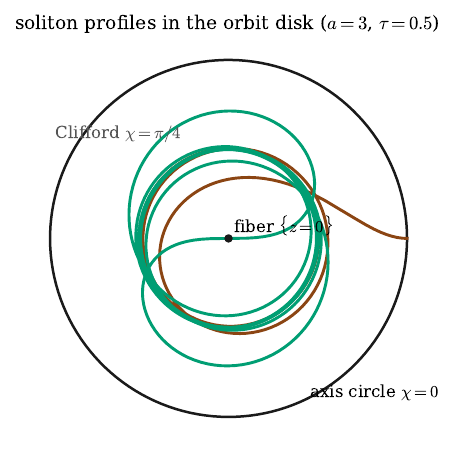}
\caption{Soliton profiles in the orbit disk $B$ ($\tau=\tfrac12$, $a=3$), drawn in
the polar chart $r=\tfrac\pi2-\chi$: the boundary circle is the axis $\{w=0\}$, the
center is the polar fiber $\{z=0\}$, and the dashed circle is the Clifford profile
$\chi=\pi/4$. Brown: axis-crossing profile; green: the two arms of the
fiber-crossing profile. All spiral onto the Clifford circle.}
\label{fig:profilesdisk}
\end{figure}

\begin{figure}[t]
\includegraphics[width=\linewidth]{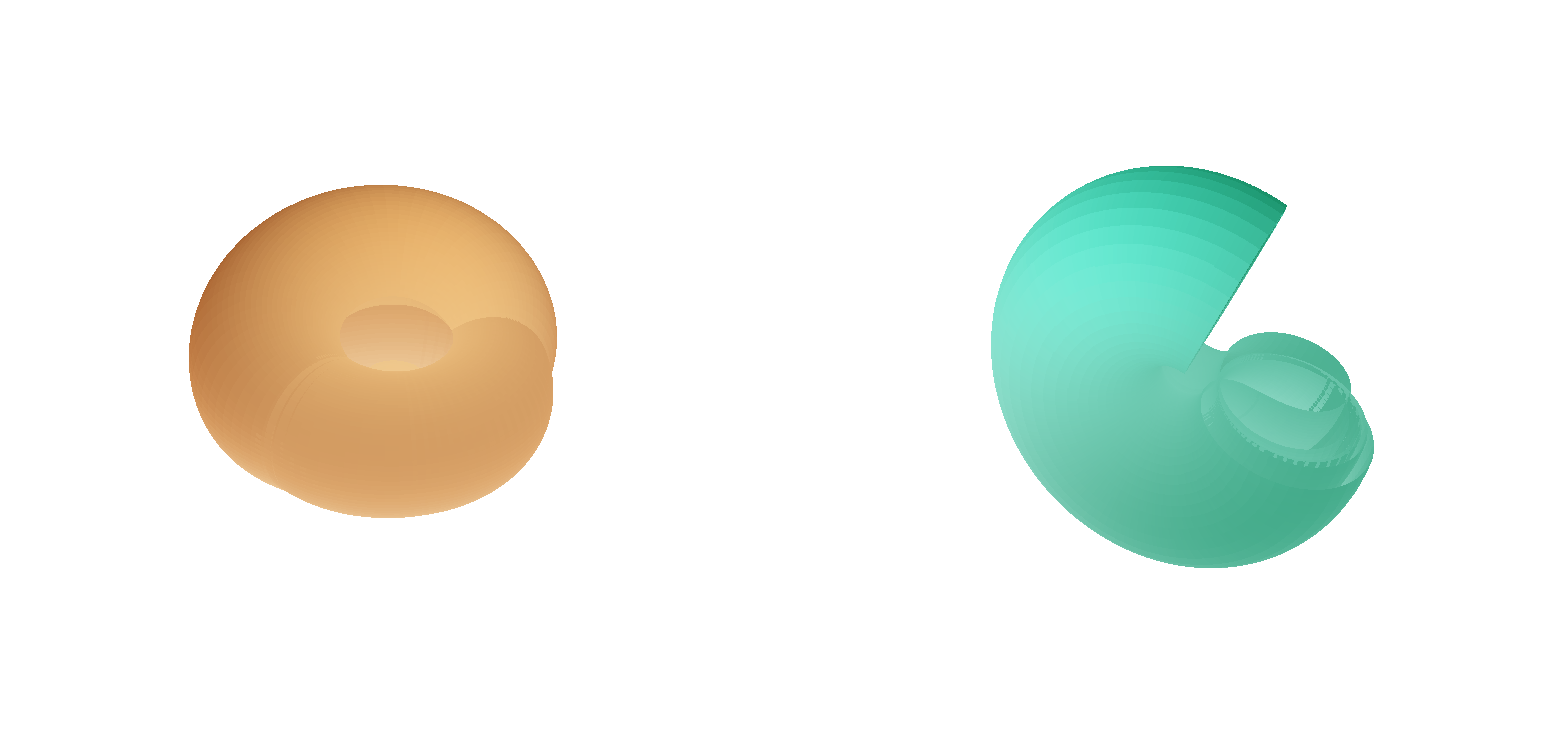}
\caption{Stereographic images (cutaway) of the two exceptional rotational solitons
for $\tau=\tfrac12$, $a=3$: the axis-crossing soliton disk (left) and the
fiber-crossing soliton cylinder (right). Both ends wrap onto the Clifford torus.}
\label{fig:rotsurfaces}
\end{figure}

(5) \emph{Axis:} a profile meeting $\{\chi=0\}$ must do so orthogonally
($\sin\beta\to0$), for otherwise $\beta'\sim-\sin\beta/\chi$ forces
$|\beta|\to\infty$, incompatible with bounded $\beta'$ elsewhere; the surface is then
smooth across the fixed point (the standard removable-singularity for equivariant
surfaces meeting the axis orthogonally). Near $\chi=0$ with $\beta(0)=0$ the
equation $\beta'=-\beta/\sigma+a\tau+O(\sigma)$ has the unique bounded solution
$\beta=\tfrac{a\tau}{2}\sigma+O(\sigma^2)$ (Theorem \ref{thm:picard} applied to the
regular singular extension across $\sigma=0$): exactly one orbit emanates from each
boundary point, and boundary points are all equivalent under the rotations
$\Phi^{1,0}_t$. Topologically the surface is a disk: the profile arc, reflected at
the axis, is a single curve whose $\mathcal R$-orbit is an immersed plane.
\emph{Polar fiber:} the profile passes through the interior point $\chi=\pi/2$ of
$B$ with a well-defined direction; all directions are equivalent under
$\Phi^{1,0}_t$, and through a point with given direction passes exactly one
solution (Theorem \ref{thm:picard}): uniqueness up to isometry. The preimage of a
profile through the polar fiber is a cylinder containing that fiber. Both exceptional
profiles were computed numerically in Figures \ref{fig:phaserot}--\ref{fig:rotsurfaces}.
Since a generic profile has its unique degenerate event of tangency type
$\sin\beta=0$ with $\chi\in(0,\pi/2)$, it avoids $\chi\in\{0,\pi/2\}$ and generates an
immersed cylinder.

(1) Existence and uniqueness of $\gamma_{\chi_0}$ for each $\chi_0\in[0,\pi/2]$ is
Theorem \ref{thm:picard} applied to \eqref{eq:rotsys} on the open cylinder (and, at
the two endpoints, to its regular singular extensions computed in item (5));
completeness of each $\gamma_{\chi_0}$ is item (4). That every profile arises this
way is item (3). Finally, $\chi_0$ is an isometry invariant of $\gamma_{\chi_0}$
under $\theta\mapsto\theta+c$ and $\theta\mapsto-\theta$ --- the isometries
normalizing the class of $a\partial_\theta$-solitons --- because both fix $\chi$
identically; hence $\chi_0\mapsto\gamma_{\chi_0}$ is injective on congruence classes,
and $\{\gamma_{\chi_0}\}_{\chi_0\in[0,\pi/2]}$ is exactly the moduli space of
nontrivial complete rotational solitons for the given $(a,\tau)$.
\end{proof}

\begin{remark}[dependence on $\tau$]\label{rem:rottau}
Unlike the Hopf class, the rotational class genuinely feels the Berger deformation:
the profile equation \eqref{eq:rotsys} and the threshold $a_*=4\sqrt{1+\tau^2}/\tau$
depend on $\tau$. As $\tau\to1$ we recover the picture in the round sphere with
$a_*=4\sqrt2$; as $\tau\to0^+$, $a_*\to\infty$, so in strongly collapsed Berger
spheres every rotational soliton end spirals. Note also that
$\mu=a\sqrt{q(\pi/4)}=a\,|\partial_\theta|_{\hat g}(\pi/4)$ is exactly the speed of
the rotation measured on the Clifford profile, which explains the form of
\eqref{eq:roteigen}.
\end{remark}

\begin{remark}[Hopf field solitons]
The Hopf field itself, $V=X_{1,1}$, has $a=1$: rotational surfaces which are
$V$-solitons (``Hopf rotators'') satisfy \eqref{eq:rotsys} with $a=1<a_*$ for every
$\tau$, hence always spiral. By contrast, for Hopf tubes $V$ acts trivially,
$\alpha=0$: the $V$-solitons among Hopf tubes are exactly the minimal ones,
$\pi^{-1}(\text{great circle})$.
\end{remark}

\section*{Statements and Declarations}

\subsection*{Funding}
The authors declare that no funds, grants, or other support were received during the
preparation of this manuscript.

\subsection*{Competing Interests}
The authors have no relevant financial or non-financial interests to disclose.

\subsection*{Data Availability}
The Python code used to produce the numerical phase portraits and stereographic
figures in this article is available from the corresponding author upon reasonable
request.

\subsection*{Author Contributions}
All authors contributed equally to the conception of the problem, the proofs, the
numerical verification, and the writing of the manuscript. All authors read and
approved the final manuscript.

\end{document}